\documentclass[12pt,reqno]{amsart}
\usepackage{graphicx}
\usepackage{amssymb,amsmath}
\usepackage{amsthm}
\usepackage{hyperref}
\usepackage{color}

\numberwithin{equation}{section}
\numberwithin{figure}{section}

\DeclareMathOperator{\R}{\mathbb{R}}

\begin{document}
		
\title[Modular Burgers' equation]{\bf Extinction of multiple shocks \\ in the modular Burgers' equation}

\author{Dmitry E. Pelinovsky}
\address[D. Pelinovsky]{Department of Mathematics and Statistics, McMaster University, Hamilton, Ontario, Canada, L8S 4K1}
\email{dmpeli@math.mcmaster.ca}

\author{Bj\"{o}rn de Rijk}
\address[B. de Rijk]{Karlsruhe Institute for Technology, Englerstra\ss e 2, 76131 Karlsruhe, Germany}
\email{bjoern.de-rijk@kit.edu}


\begin{abstract}
	We consider multiple shock waves in the Burgers' equation with a modular advection term. It was previously shown that the modular Burgers' equation admits a traveling viscous shock with a single interface, which is stable against smooth and exponentially localized perturbations. In contrast, we suggest in the present work with the help of energy estimates and numerical simulations that the evolution of shock waves with multiple interfaces leads to finite-time coalescence of two consecutive interfaces. We formulate a precise scaling law of the finite-time extinction supported by the interface equations and by numerical simulations. 
\end{abstract}

\date{\today}
\maketitle


\section{Introduction}

The present work addresses multiple shock waves in the modular Burgers' equation 
\begin{equation}
\label{Burgers}
\frac{\partial u}{\partial t} = \frac{\partial |u|}{\partial x} + \frac{\partial^2 u}{\partial x^2},
\end{equation}
which is different from the classical Burgers' equation by the modular advection 
term. Equation~\eqref{Burgers} has been used as a model to describe inelastic dynamics 
of particles with piecewise interaction potentials~\cite{Hedberg1,RH18}. 
Generalizations of this model with additional terms were also discussed in~\cite{RH16,Hedberg2,Hedberg3}. 

Some preliminary results were obtained for the modular Burgers' equation~(\ref{Burgers}) both analytically and numerically. Traveling wave solutions were constructed in~\cite{Rudenko1,Rudenko3} by matching solutions at the interfaces where the modular nonlinearity jumps. Collisions of compactly supported pulses and dynamics near a viscous shock were studied in~\cite{Hedberg1} by using qualitative approximations. Numerical approximations of time-dependent solutions of the modular Burgers' equation~(\ref{Burgers}) were constructed in~\cite{Radostin2} with the aid of Fourier sine series. 

A traveling viscous shock of the form $u(t,x) = U_c(x-ct)$ is available in the exact analytical form:
\begin{equation}
\label{shock}
U_c(\xi) = \left\{ \begin{array}{ll} U_+ (1 - e^{(1+c)(\xi_0 - \xi)}), \quad & \xi > \xi_0, \\
U_- (1 - e^{(1-c)(\xi-\xi_0)}), \quad & \xi \leq \xi_0, \end{array}
\right.
\end{equation}
where $\lim\limits_{\xi \to \pm \infty} U_c(\xi) = U_{\pm}$ satisfy $U_- < 0 < U_+$ and the speed $c$ is uniquely selected at 
\begin{equation}
\label{speed}
c = -\frac{U_+ + U_-}{U_+ - U_-}.
\end{equation}
The viscous shock~(\ref{shock}) possesses a single interface at $\xi_0 \in \R$, where $U_c(\xi_0) = 0$, such that $U_c$ is a continuously differentiable function with a piecewise continuous second derivative, whose only discontinuity arises at the interface and is described by the jump condition:
\begin{equation}
\label{interface-ode}
[U_c'']^+_-(\xi_0) := U_c''(\xi_0+0) - U_c''(\xi_0-0) = - 2 |U_c'(\xi_0)|.
\end{equation}

Asymptotic stability of the traveling shock~(\ref{shock}) against smooth, exponentially localized perturbations was established in~\cite{LPP21}. It was shown that the evolution of such perturbations is well defined on both sides of the interface and the perturbations decay in time. A finite-difference numerical method, which couples the nonlinear dynamics at the interfaces to the linear advection-diffusion dynamics on both sides of the interface, was implemented in~\cite{LPP21} to corroborate the stability analysis of the viscous shock. 

The purpose of the present work is to study viscous shocks in the modular Burgers' equation~(\ref{Burgers}) with multiple interfaces. We show with the aid of energy estimates that compact regions between two consecutive interfaces shrink in time and no new compact regions can be formed dynamically. In particular, this yields that no traveling viscous shocks with multiple interfaces can exist in the modular Burgers' equation. Moreover, for odd initial data with three symmetric interfaces we establish that the interfaces coalesce in finite time to a single interface. We complement our analysis with finite-difference numerical simulations. Postprocessing data analysis suggests a precise scaling law of the finite-time extinction which agrees with the interface equations.

We note that, although the finite-difference method is rather elementary, it allows us to capture the main feature of the dynamics of 
the modular Burgers' equation~(\ref{Burgers}), where the linear equations 
between interfaces are coupled together by the nonlinear 
interface equations. It is unclear how else the numerical modeling of the time evolution could be performed due to the singular contribution of 
the modular nonlinearity (without replacing it by a smooth approximation). 

Before closing the introduction, we mention some contemporary work 
on other related problems. A diffusion equation with a piecewise defined nonlinearity, namely, the KPP model with the cutoff reaction rate, was studied in~\cite{KPP3,KPP4}, where matched asymptotic expansions in the dynamically moving coordinate frame have been used to establish both the existence and asymptotic stability of traveling viscous shocks. Metastable $N$-waves of the classical Burgers' equation were studied in~\cite{B2,B3} by employing dynamical systems methods. 

The ultimate goal of our studies is to understand the dynamics of the logarithmic Burgers' equations~\cite{J20}, which commonly arises in the modeling 
of granular chains in viscous systems. The logarithmic nonlinearity is
more singular than the modular nonlinearity in~\eqref{Burgers}, hence it presents further challenges in the analysis of (traveling) viscous shocks. We remark that, compared to the logarithmic Burgers' equation, the logarithmic diffusion equation has been well-studied~\cite{Carles}.

The paper is organized as follows. Section~\ref{sec-3} contains general 
energy estimates for the modular Burgers' equation. Section~\ref{sec-new}
addresses the finite-time extinction of shocks for odd initial data with three symmetric interfaces, both analytically and numerically. 
Section~\ref{sec-6} concludes the paper with a discussion of open problems.

\section{Energy estimates}
\label{sec-3}

Here we use energy estimates to show that a compact region between two consecutive interfaces shrinks and eventually disappears in the time evolution of the modular Burgers' equation~(\ref{Burgers}). We take $T > 0$ and consider a continuously differentiable solution $u(t,x) \colon (0,T) \times \R \to \R$ to the modular Burgers' equation~\eqref{Burgers}, whose second derivative $u_{xx}$ is piecewise continuous with discontinuities arising only at the interfaces. 

We consider two consecutive interfaces $-\infty < \xi_1(t) < \xi_2(t) < \infty$, so that $u(t,\xi_1(t)) = 0$ and $u(t,\xi_2(t)) = 0$ for $t \in (0,T)$. Without loss of generality, we assume $u(t,x) > 0$ for $\xi_1(t) < x < \xi_2(t)$. All in all, this yields the following linear evolutionary boundary-value problem:
\begin{equation}
\label{lin-Bur}
\left\{ 
\begin{array}{ll}
u_t = u_x + u_{xx}, \quad \xi_1(t) < x < \xi_2(t), & \quad 0 < t < T, \\
u(t,\xi_1(t)) = 0, & \quad 0 < t < T, \\
u(t,\xi_2(t)) = 0, & \quad 0 < t < T.
\end{array}
\right.
\end{equation}
The linear problem~(\ref{lin-Bur}) is not closed as we need to find the evolution of $\xi_{1,2}(t)$ from the boundary conditions at $x = \xi_{1,2}(t)$ and the evolutionary boundary-value problems 
satisfied by $u(t,x)$ for $x < \xi_1(t)$ and for $x > \xi_2(t)$. At each interface $x = \xi_{1,2}(t)$, two additional boundary conditions are needed. These two conditions are given by the continuity of the derivative $u_x(t,x)$ across $x = \xi_{1,2}(t)$ and by a jump condition for $u_{xx}(t,x)$, which read
\begin{align}
\label{interface-eq}
\begin{split}
[u_{x}]^+_-(t,\xi_{1,2}(t)) &:= u_{x}(t,\xi_{1,2}(t)+0) - u_{x}(t,\xi_{1,2}(t)-0) = 0,\\
[u_{xx}]^+_-(t,\xi_{1,2}(t)) &:= u_{xx}(t,\xi_{1,2}(t)+0) - u_{xx}(t,\xi_{1,2}(t)-0) = -2|u_x(t,\xi_{1,2}(t))|,
\end{split}
\end{align}
for $0 < t < T$. We note that the jump condition for $u_{xx}$ is equivalent to the continuity of the temporal derivative $u_t$ across the interfaces. 
We derive energy estimates from the linear boundary-value problem~(\ref{lin-Bur}) by ignoring the global information from other boundary conditions~(\ref{interface-eq}). Consequently, the time evolution of $\xi_{1,2}(t)$ is not relevant for our energy estimates. 

Integrating~(\ref{lin-Bur}) on $[\xi_1(t),\xi_2(t)]$ yields 
\begin{equation}
\label{reduction-1}
\frac{d}{dt} \int_{\xi_1(t)}^{\xi_2(t)} u(t,x) dx = u_x(t,\xi_2(t)) - u_x(t,\xi_1(t)) \leq 0,
\end{equation}
where we have used the inequalities $u_x(t,\xi_2(t)) \leq 0$ and 
$u_x(t,\xi_1(t)) \geq 0$, which follow from the fact that
$u(t,x) > 0$ for $\xi_1(t) < x < \xi_2(t)$. Hence, 
the positive mass $\smash{\int_{\xi_1(t)}^{\xi_2(t)} u(t,x) dx}$ is monotonically 
decreasing as a function of $t$ as long as the slopes 
at the end points of the compact region are nonzero.

Integrating~(\ref{lin-Bur}) multiplied by $u$ on $[\xi_1(t),\xi_2(t)]$ yields 
\begin{equation}
\label{reduction-2}
\frac{d}{dt} \int_{\xi_1(t)}^{\xi_2(t)} u^2(t,x) dx = - 2 \int_{\xi_1(t)}^{\xi_2(t)} u_x^2(t,x) dx \leq 0.
\end{equation}
Hence, the positive energy $\smash{\int_{\xi_1(t)}^{\xi_2(t)} u^2(t,x) dx}$ is monotonically decreasing as a function of $t$ as long as $\xi_1(t) < \xi_2(t)$. In particular, identity~\eqref{reduction-2} shows that no traveling viscous shocks with multiple interfaces can exist in the modular Burgers' equation~\eqref{Burgers}, as for such solutions the positive energies $\smash{\int_{\xi_1(t)}^{\xi_2(t)} u^2(t,x) dx}$ and $\smash{\int_{\xi_1(t)}^{\xi_2(t)} u_x^2(t,x) dx}$ between two consecutive interfaces must stay constant in time.

We remark that energy estimates involving spatial derivatives of $u(t,x)$ cannot be derived from the linear boundary-value problem~(\ref{lin-Bur}), because of the lack of information on the spatial derivatives of $u(t,x)$ at $x = \xi_{1,2}(t)$. 

The two estimates~(\ref{reduction-1}) and~(\ref{reduction-2}) suggest that no new compact regions may be formed dynamically in time since the mass and energy of the compact region with positive $u(t,x)$ cannot increase from zero to positive values. However, the argument does not clarify if the mass and energy extinguish in finite or infinite times or if the two interface $\xi_{1,2}(t)$ coalesce when the mass and energy vanish. In the next section we will answer these questions for the special case of odd shock waves and corroborate our analysis by numerical experiments. 

\section{Odd initial data with three symmetric interfaces}
\label{sec-new}

\begin{figure}[htbp] 
	\centering
	\includegraphics[width=3in, height = 2in]{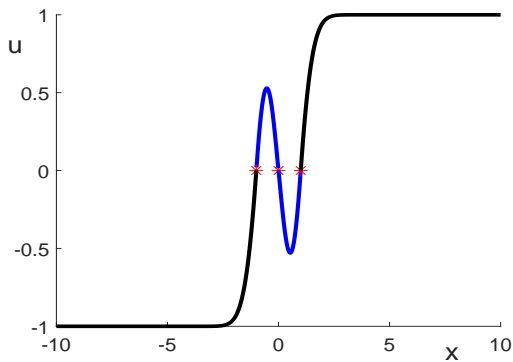}
	\caption{Odd initial data with three symmetric interfaces shown by red stars. The depicted initial data on $(0,\infty)$ is given by~(\ref{initial-data}) with $\alpha = 1$.}
	\label{fig-data}
\end{figure}

Here we consider the simplest problem for shock waves with multiple interfaces. Since the modular Burgers' equation~(\ref{Burgers}) preserves odd functions in the time evolution, we restrict solutions to the class of odd functions $u(t,-x) = -u(t,x)$ closed on $(0,\infty)$ subject to Dirichlet condition at $x = 0$ and the normalized boundary condition $u(t,x) \to 1$ as $x \to +\infty$. We will assume that there exists a single interface at $x = \xi(t) \in (0,\infty)$. Due to the oddness condition, the multiple shock wave consists of three symmetric interfaces at $x = -\xi(t)$, $x = 0$ and $x = \xi(t)$, cf.~Figure~\ref{fig-data}.

The mathematical formulation of the evolutionary boundary-value problem is given by 
\begin{equation}
\label{num-Burg}
\left\{ 
\begin{array}{lll}
u_t = -u_x + u_{xx}, & \quad u(t,x) < 0, & \quad 0 < x < \xi(t), \\
u_t = u_x + u_{xx}, & \quad u(t,x) > 0,  & \quad \xi(t) < x < \infty, \\
u(t,0) = 0, & \quad u(t,\xi(t)) = 0, & \quad \lim\limits_{x \to +\infty} u(t,x) = 1,
\end{array}
\right.
\end{equation}
where $u(t,x) : (0,T) \times \R \to \R$ is odd in $x$ and continuously differentiable with piecewise continuous second derivative $u_{xx}$ having discontinuities only at the interfaces $x = 0,\pm \xi(t)$. The supplementary conditions at the interfaces read 
\begin{equation}
\label{num-interface}
\begin{split}
u_{xx}(t,+0) &= u_x(t,0), \\
[u_{x}]^+_-(t,\xi(t)) &= 0, \\
[u_{xx}]^+_-(t,\xi(t)) &= -2 u_x(t,\xi(t)).
\end{split}
\end{equation}
The first interface condition in~(\ref{num-interface}) is consistent with the Dirichlet condition $u(t,0) = 0$ and the evolution system in~(\ref{num-Burg}). The last interface condition in~(\ref{num-interface}) in combination with the Dirichlet condition 
$u(t,\xi(t)) = 0$ and the evolution system in~(\ref{num-Burg}) can be rewritten as 
the differential equation
\begin{equation}
\label{dynamics-interface}
\xi'(t) = -1 - \frac{u_{xx}(t,\xi(t)+0)}{u_x(t,\xi(t))} = 
+ 1 - \frac{u_{xx}(t,\xi(t)-0)}{u_x(t,\xi(t))}.
\end{equation}
The interface equation (\ref{dynamics-interface}) holds as long as $u_x(t,\xi(t)) > 0$.

\subsection{Finite-time coalescence of interfaces}
\label{sec-4.1}

We establish finite-time coalescence of interfaces for solutions $u(t,x) \colon (0,T) \times \R \to \R$ of the boundary-value problem~\eqref{num-Burg}-\eqref{num-interface} by introducing the new variable $z(t,x) = 1 - u(t,x)$, which measures the difference between the solution $u(t,x)$ and its asymptotic value $\smash{\displaystyle\lim_{x \to +\infty} u(t,x)} = 1$. Clearly, $z$ satisfies the equation
\begin{align*}\frac{\partial z}{\partial t} = \frac{\partial |1-z|}{\partial x} + \frac{\partial^2 z}{\partial x^2}.\end{align*}
We aim to derive a differential inequality for the mass 
\begin{align*} M(t) = \int_0^\infty z(t,x) dx.\end{align*}
In order to assure that the mass is well-defined and positive, we require $z(t,\cdot)$ to be integrable and nonnegative on $(0,\infty)$ for each $t \in (0,T)$. By a standard application of  maximum and comparison principles for advection-diffusion equations, this is the case if the initial condition $z(0,x) = 1 - u(0,x)$ is nonnegative and integrable on $(0,\infty)$. 

The advantage of working with the mass $M(t)$ over $\smash{\int_0^{\xi(t)} u(t,x) dx}$, as in Section~\ref{sec-3}, is that information on \emph{both} sides of the interface at $x = \xi(t)$ is taken into account. That is, the estimate~\eqref{reduction-1} only relies on the \emph{linear} dynamics~\eqref{lin-Bur} between interfaces, whereas the modular nonlinearity can only be captured by considering both sides of an interface.

Let us proceed with deriving a differential inequality for the mass $M(t)$. We take $t \in (0,T)$ and assume no coalescence of interfaces has occurred on $[0,t]$, so that it holds $\xi(s) > 0$ for all $s \in [0,t]$. First, we note that 
$$
z(s,\xi(s)) = z(s,0) = 1, \qquad -z_x(s,0) = u_x(s,0) \leq 0,
$$ 
and 
$$
\smash{\displaystyle \lim_{x \to +\infty} z(s,x)} = \smash{\displaystyle \lim_{x \to +\infty} z_x(s,x)} = 0
$$ 
hold for all $s \in (0,t)$. Hence, with the aid of the Leibniz rule, we obtain
\begin{align*} \frac{d}{ds} M(s) &= \lim_{x \to +\infty} \frac{d}{ds} \left(\int_0^{\xi(s)} z(s,y) dy + \int_{\xi(s)}^x z(s,y) dy\right)\\
&= \lim_{x \to +\infty} \left(z(s,\xi(s))\xi'(s) + \int_0^{\xi(s)} \left(z_{xx}(s,y) - z_x(s,y)\right) d y\right.\\ 
&\qquad \qquad \qquad \left. \, - z(s,\xi(s))\xi'(s) + \int_{\xi(s)}^x \left(z_{xx}(s,y) + z_x(s,y)\right) dy\right)\\
&= \lim_{x \to +\infty} \left(z_{x}(s,\xi(s)) - 1 - \left(z_x(s,0) - 1\right) + z_x(s,x) + z(s,x) - \left(z_x(s,\xi(s)) + 1 \right)\right)\\
&= -1-z_x(s,0) \leq -1
\end{align*}
for $s \in (0,t)$, where we remark that the interchange of limit and derivative is justified by uniform convergence of the relevant differential quotient. Upon integrating the above differential inequality for the mass $M(s)$ for $s \in [0,t]$, we obtain $M(t) \leq M(0) - t$.

To finish the argument, we note that $z(t,x) \geq 1$ for $x \in (0,\xi(t))$ and $z(t,x) \geq 0$ for $x \in (\xi(t),\infty)$, so that we arrive at
\begin{align}
\label{time-estimate} 
0 < \xi(t) \leq \int_0^{\xi(t)} z(t,x) dx \leq M(t) \leq M(0) - t,
\end{align}
which implies that there exists $t_0 \in (0,M(0))$ such that $\xi(t) \to 0$ as $t \to t_0$. Hence we have established finite-time coalescence of interfaces for the simplest odd multiple shock waves. Moreover, our method provides an upper bound for the time of coalescence, which is given by the integral
\begin{align*} M(0) = \int_0^\infty \left(1 - u_0(x)\right) dx,\end{align*} 
where $u_0(x) = u(0,x)$ is the initial condition of the solution $u(t,x)$ to the evolutionary boundary-value problem~\eqref{num-Burg}-\eqref{num-interface}.

\subsection{Finite-difference method}
\label{sec-4}

Next, we set up the framework for the numerical experiments, which rely on a finite-difference method. To implement the method we assume that $\xi(t) > 0$ and work with the rescaled spatial coordinate $y := x/\xi(t)$. This transformation scales 
the domain of the boundary-value problem~(\ref{num-Burg}) to the time-independent regions $(0,1)$ and $(1,\infty)$. Thus, abusing notation, 
we rewrite the evolutionary boundary-value problem for $u = u(t,y)$ as:
\begin{equation}
\label{num-Burg-new}
\left\{ 
\begin{array}{lll}
u_t = \xi^{-1} (\xi' y - 1) u_y + \xi^{-2} u_{yy},  & \quad u(t,y) < 0, & \quad 0 < y < 1, \\
u_t = \xi^{-1} (\xi' y + 1) u_y + \xi^{-2} u_{yy},  & \quad u(t,y) > 0, & \quad 1 < y < \infty, \\
u(t,0) = 0, & \quad u(t,1) = 0, & \quad \lim\limits_{y \to +\infty} u(t,y) = 1,
\end{array}
\right.
\end{equation}
whereas the interface equation~(\ref{dynamics-interface}) transforms into 
\begin{equation}
\label{dynamics-interface-new}
\xi'(t) = -1 - \frac{u_{yy}(t,1+0)}{\xi(t) u_y(t,1)} = 
+ 1 - \frac{u_{yy}(t,1-0)}{\xi(t) u_y(t,1)}.
\end{equation}

By using an equally spaced grid with the step size $h$ on $[0,1]$ and $[1,L]$ for sufficiently large $L$, we replace the first and second spatial derivatives in~(\ref{num-Burg-new}) by the central differences. We can do this for every interior point of the grid 
since  there are no evolution equations at $y = 0$ and $y = 1$ due to the Dirichlet conditions. The Neumann condition $u_y(t,L) = 0$ is used at $y = L$. It remains to derive a discretization of the interface condition~(\ref{dynamics-interface-new}).

To couple the solutions on $[0,1]$ and $[1,L]$, we use the central difference approximation of the first and second spatial derivatives at $y = 1$ in~(\ref{dynamics-interface-new}). This can only be done if additional grid points are added to the left and to the right of the interface point $y = 1$. In other words, we augment $\{ u_k \}_{k=0}^{k=N}$ for $y_k = hk$ with $v_{N+1}$ for $y_{N+1} = 1+h$ and $\{ u_k \}_{k=N}^{k=M}$ for $y_k = hk$ with $v_{N-1}$ for $y_{N-1} = 1 - h$, where $h = \frac{1}{N} = \frac{L}{M}$ and $u_0 = u_N = 0$ due to the Dirichlet conditions at $y = 0$ and $y = 1$. For the Neumann condition at $y = L$, we use an additional grid point at $y_{M+1} = L+h$ with $u_{M+1} = u_{M-1}$.

The continuity of $u_y(t,y)$ and the jump of $u_{yy}(t,y)$ across $y = 1$ 
are expressed in the central difference approximation by the linear equations
\begin{align*}
\frac{v_{N+1} - u_{N-1}}{2h} &= \frac{u_{N+1} - v_{N-1}}{2h}, \\
\frac{u_{N+1} + v_{N-1}}{h^2} - \frac{v_{N+1} + u_{N-1}}{h^2} &= - 2\xi \frac{v_{N+1} - u_{N-1}}{2h}.
\end{align*}
These linear equations admit a unique solution for the additional variables $v_{N+1}$ and $v_{N-1}$ given by 
\begin{align*}
v_{N+1} = \frac{2 u_{N+1} - h \xi u_{N-1}}{2 - h \xi}, \quad 
v_{N-1} = \frac{2 u_{N-1} - h \xi u_{N+1}}{2 - h \xi},
\end{align*}
where we assume that $h$ is chosen so small that $h \xi(t) < 2$. Substituting these solutions into the central difference approximation of the interface equation~(\ref{dynamics-interface-new}) yields the approximation 
\begin{equation}
\label{dynamics-interface-num}
\xi'(t) = -\frac{(2 - h \xi) (u_{N+1} + u_{N-1})}{h \xi (u_{N+1} - u_{N-1})}.
\end{equation}

The time evolution of the linear system~(\ref{num-Burg-new})  was approximated by the implicit Crank-Nicolson method based on the trapezoidal rule of numerical integration. The Crank-Nicolson method is unconditionally stable for the linear advection-diffusion equations. 
However, the stability of iterations was affected by the approximation 
(\ref{dynamics-interface-num}) since $\xi(t)$ and $\xi'(t)$ were used in the evolutionary system~(\ref{num-Burg-new}) explicitly based on the predictor-corrector pair (with $\xi(t)$ obtained from $\xi'(t)$ by using Heun's method). 

It remains to provide initial data $u_0(x) := u(0,x)$, which are consistent with the interface conditions~(\ref{dynamics-interface-new}). Without loss of generality, we assume $\xi(0) = 1$ so that $y = x$ at $t = 0$. For sufficiently fast convergence of $u_0(x)$ towards $1$ as $x \to +\infty$, we consider a Gaussian function on $(1,\infty)$ concatenated with a quartic polynomial on $(0,1)$:
\begin{align}
\label{initial-data}
u_0(x) = \left\{ \begin{array}{ll} x (1-x) (ax^2 + bx + c), & \quad 0 < x < 1, \\
1 - e^{-\alpha (x^2 - 1)}, &\quad 1 < x < \infty, \end{array} \right.
\end{align}
so that the boundary conditions $u_0(0) = u_0(1) = 0$ and $\lim\limits_{x \to +\infty} u_0(x) = 1$ are satisfied. Parameters $a$, $b$, and $c$ can then be found uniquely in terms of $\alpha$ by using the interface conditions~(\ref{num-interface}). The condition $u_0''(0) = u_0'(0)$ yields $2b = 3c$. 
The condition $u_0'(1+0) = u_0'(1-0)$ yields $a + b + c = -2\alpha$. 
Finally, the condition $u_0''(1+0) - u_0''(1-0) = - 2u_0'(1)$ yields 
$2a + b = 2 \alpha^2 - \alpha$. Solving all three conditions, 
we obtain 
\begin{align*}
a = \frac{\alpha (10 \alpha+ 1)}{7}, \quad 
b = -\frac{3 \alpha (2 \alpha + 3)}{7}, \quad 
c = -\frac{2 \alpha (2 \alpha+ 3)}{7},
\end{align*}
which completes the construction of the initial condition $u_0(x)$ for arbitrary $\alpha > 0$. Since $\xi'(0) = 2 (\alpha - 1)$, the interface expands initially if $\alpha > 1$ and contracts initially if $\alpha < 1$. 

\subsection{Outcomes of numerical simulations with the initial data~(\ref{initial-data})}

We performed iterations on the domain $[0,L]$ with $L = 10$, discretized with the step size $h = 0.02$. The time step was selected to be $\tau = 0.0001$ in order to obtain better 
accuracy in the evolution of the interface $\xi(t)$ within the 
finite-difference approximation~(\ref{dynamics-interface-num}). Nevertheless, the accuracy was decreasing when $\xi(t)$ and $u_y(t,1)$ were getting smaller and iterations eventually broke up and stopped 
before $\xi(t)$ could reach $0$. This was partly related to the fact that 
the Neumann condition at the end point $L = 10$ was preserving the initial value $u_0(L) \approx 1$ for a while, after which the value of $u(t,L)$ started to decrease during the extinction stage.

\begin{figure}[htbp] 
	\centering
	\includegraphics[width=3in, height = 2in]{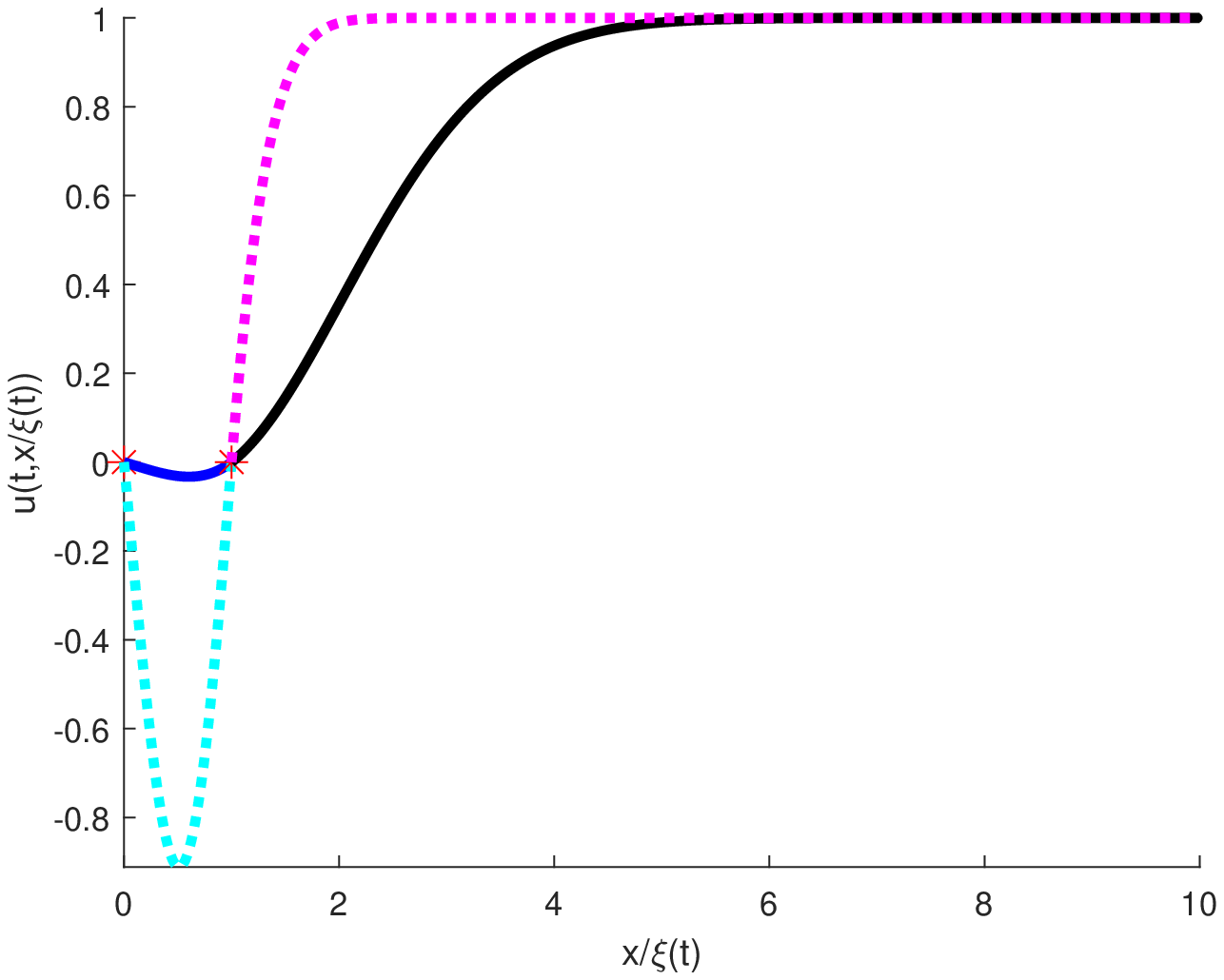}
	\includegraphics[width=3in, height = 2in]{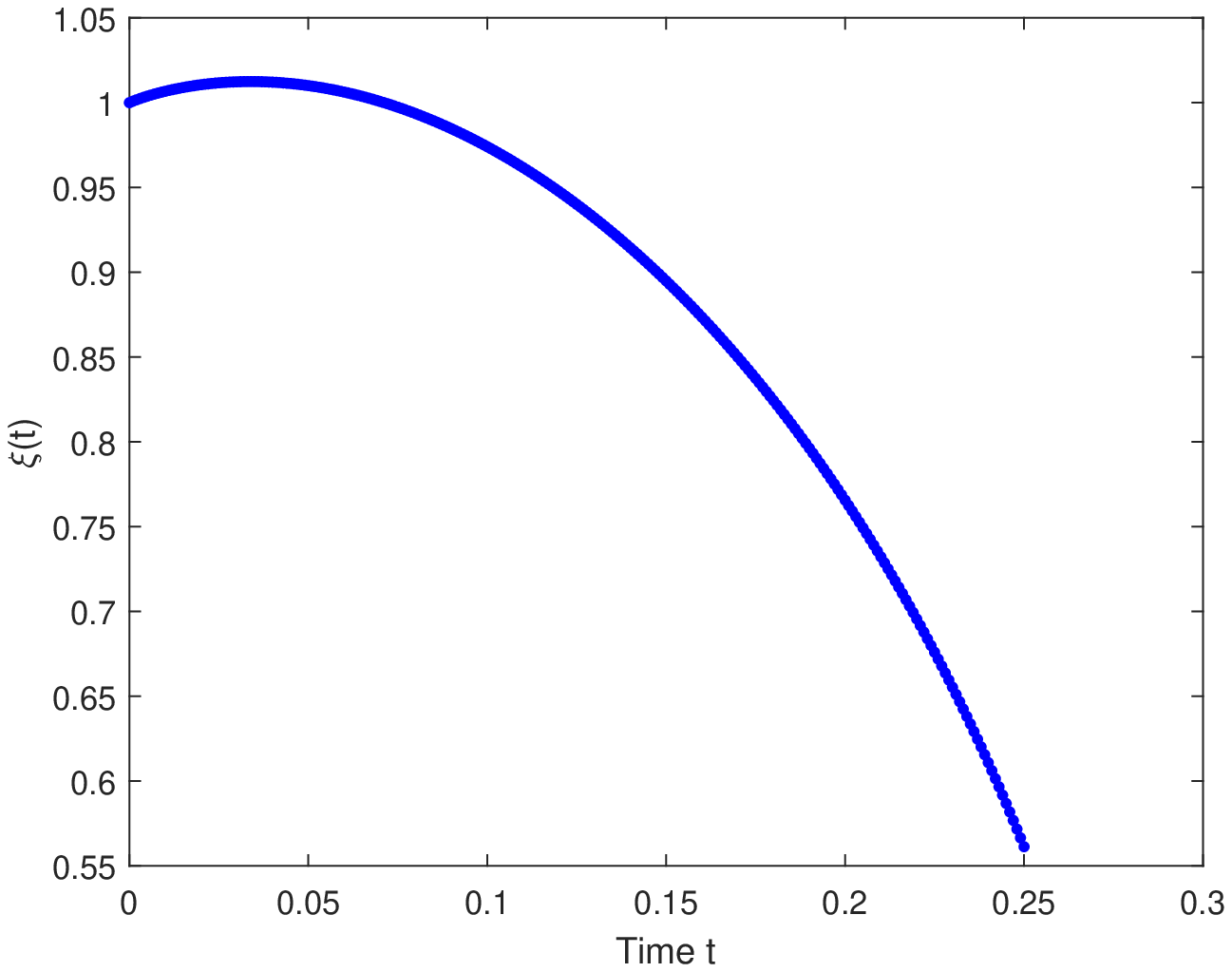} \\
	\includegraphics[width=3in, height = 2in]{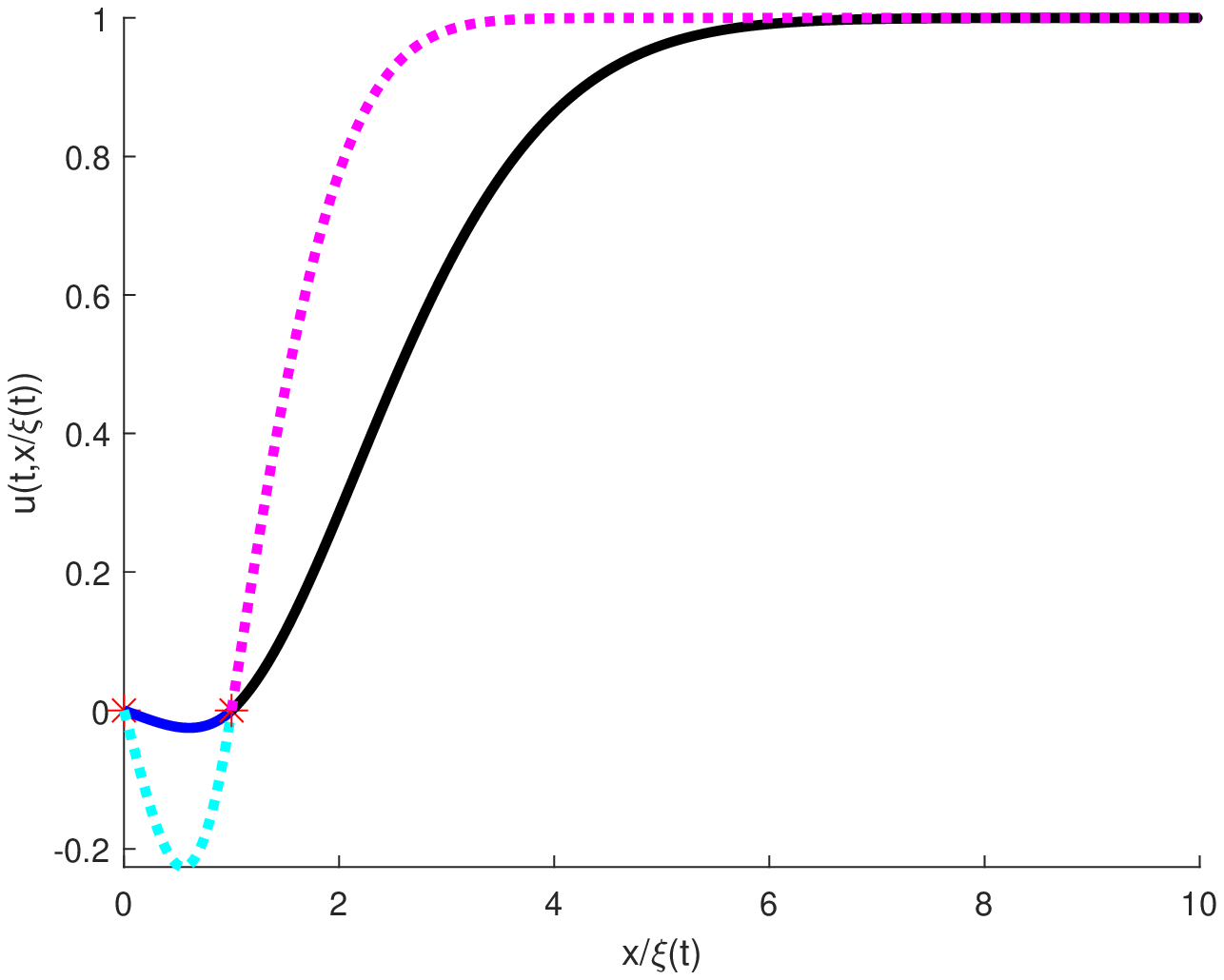}
	\includegraphics[width=3in, height = 2in]{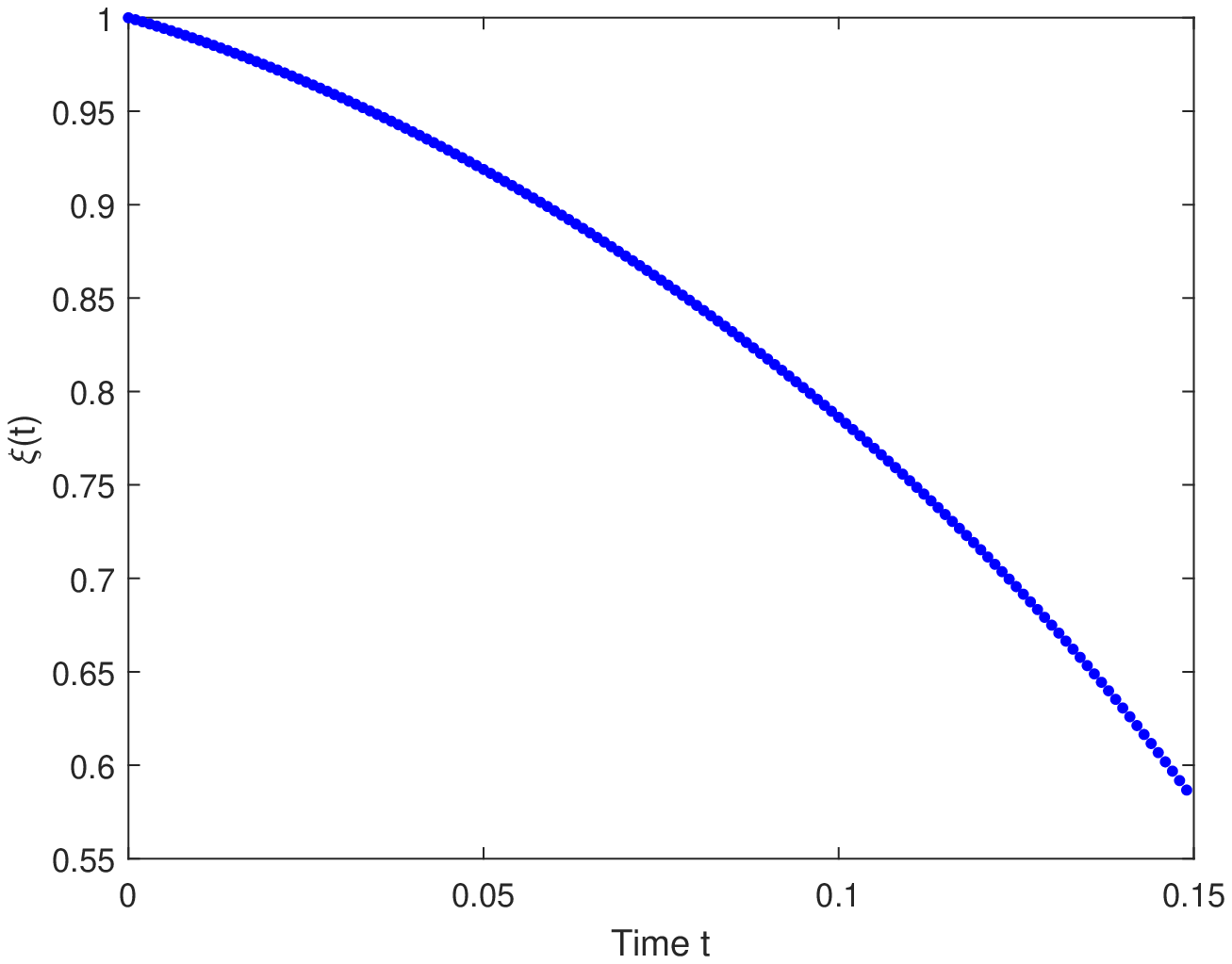} 
	\caption{Evolution of the boundary-value problem~(\ref{num-Burg-new}) for the initial data~(\ref{initial-data}) for $\alpha = 1.5$ (top) and $\alpha = 0.5$ (bottom). Left: $u(t,y)$ versus $y > 0$ for $t = 0$ (dashed) and $t = t_1$ (solid). Right: evolution of $\xi(t)$ versus $t$.}
	\label{fig-1}
\end{figure}

Figure~\ref{fig-1} depicts the outcomes of the numerical simulations with the initial data~(\ref{initial-data}) for $\alpha = 1.5$ (top) and $\alpha = 0.5$ (bottom). 
The left panels show the profile $u(t,y)$ for $y > 0$ and two values of time: $t = 0$ (dashed line) and $t = t_1$ (solid line), where $t_1 = 0.25$ for $\alpha = 1.5$ and $t_1 = 0.15$ for $\alpha = 0.5$. The right panels show the numerically computed evolution of the interface $\xi(t)$ versus $t$. 

First, we observe that the evolution of $\xi(t)$ is non-monotone for $\alpha = 1.5$ and monotone for $\alpha = 0.5$. This is in agreement with $\xi'(0) = 2(\alpha - 1)$ computed from (\ref{initial-data}). Moreover, performing computations for longer times with these and other values of $\alpha$ suggests that in all cases there exists an extinction time $t_0 \in (0,\infty)$ such that 
$$
\xi(t) \to 0, \quad u_x(t,\xi(t)) \to 0, \quad 
u_{xx}(t,\xi(t)\pm 0) \to 0, \quad {\rm as} \;\; t \to t_0,
$$ 
where the spatial derivatives were computed in the original variable $x$ by using the chain rule and the numerical approximations:
\begin{equation}
\label{num-deriv}
u_x(t,\xi(t)) = \frac{u_{N+1}-u_{N-1}}{h \xi(t) (2 - h \xi(t))}, \quad 
u_{xx}(t,\xi(t)-0) = 2\frac{u_{N+1}+u_{N-1}(1-h \xi(t))}{h^2 \xi^2(t) (2 - h \xi(t))}.
\end{equation}

Figure~\ref{fig-1} suggests that the extinction time $t_0$ of the interfaces is actually much smaller than the upper bound $T(\alpha) = M(0)$ computed 
from (\ref{time-estimate}), or explicitly,
\begin{align*} 
T(\alpha) = \int_0^\infty \left(1-u_0(x)\right) dx = \frac{\sqrt{\pi } e^\alpha \text{erfc}\left(\sqrt{\alpha}\right)}{2 \sqrt{\alpha}}+\frac{2 \alpha^2}{21}+\frac{17\alpha}{70}+1.
\end{align*}
Indeed, we find $T(1.5) \approx 1.84859$ and $T(0.5) \approx 1.80092$. Hence, the upper bound of the extinction time derived from (\ref{time-estimate}) is not sharp.

\subsection{Scaling laws describing the finite-time extinction}
\label{sec-5}

We claim based on the postprocessing data analysis that 
the following scaling law of extinction holds as $t \to t_0$:
\begin{equation}
\label{scaling-law}
\xi(t) \sim \sqrt{t_0-t}, \quad 
u_x(t,\xi(t)) \sim (t_0-t), \quad 
u_{xx}(t,\xi(t)-0) \sim \sqrt{t_0-t}.
\end{equation}
This scaling law is in agreement with the interface equation~(\ref{dynamics-interface}), which suggests that $\xi'(t)$ diverges as $t \to t_0$:
$$
\xi'(t) \sim -\frac{1}{\sqrt{t_0-t}}.
$$
For postprocessing data analysis, we use linear regression in the log-log variable, i.e. 
\begin{equation}
\label{lin-regression}
\log(\xi(t)) \quad \mbox{\rm versus} \quad c_1 \log(t_0-t) + c_2,
\end{equation}
where the coefficient $c_1$ determines the power of the scaling law~(\ref{scaling-law}). The only obstacle with this method is that the value of $t_0$ is unknown 
and cannot be approximated well because the iterations break down when $\xi(t)$ becomes too small (in our simulations smaller than $0.3$). 

\begin{figure}[htbp] 
	\centering
	\includegraphics[width=3in, height = 2in]{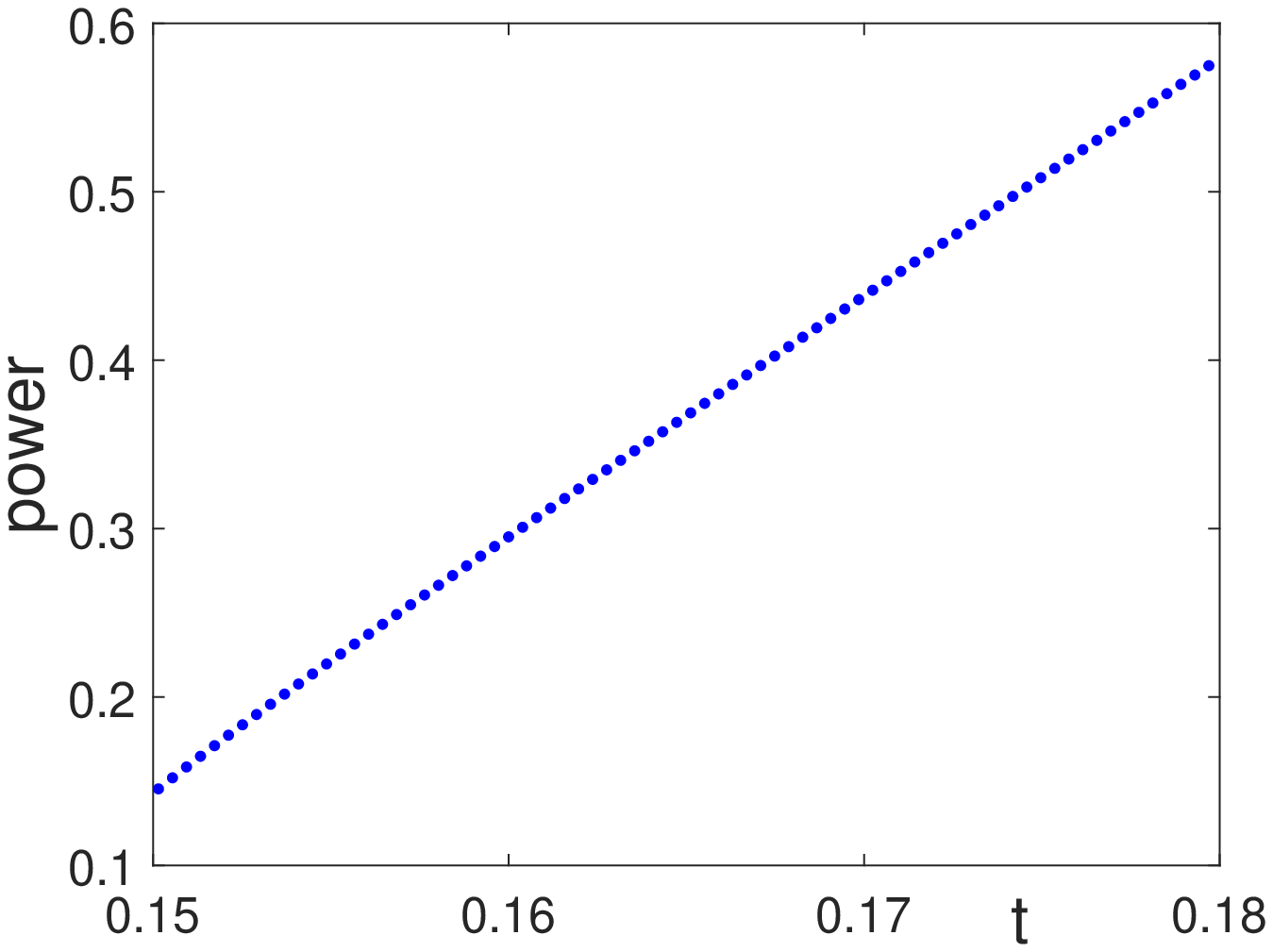}
	\includegraphics[width=3in, height = 2in]{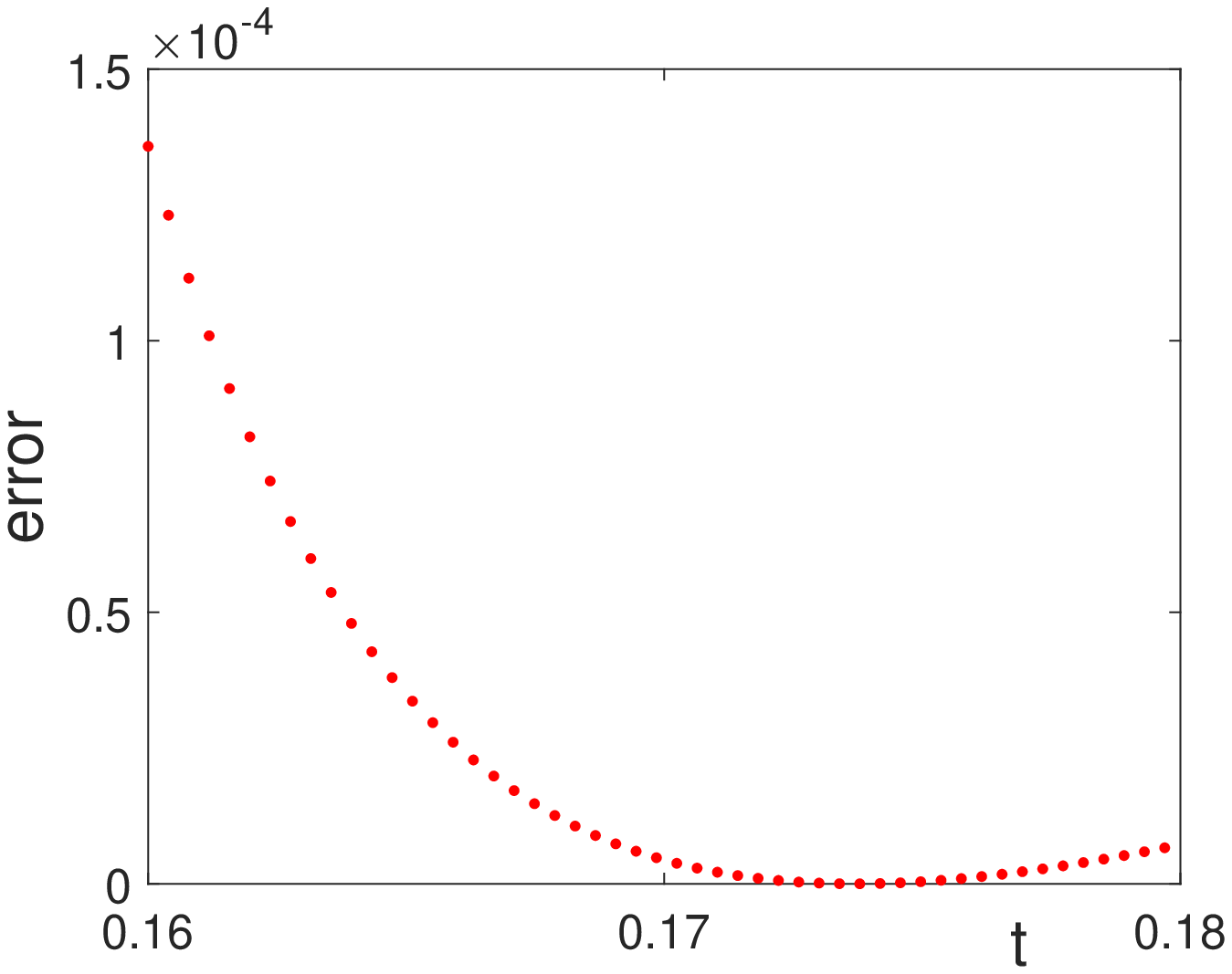} 
	\caption{Power of the linear regression (left) and the approximation error (right) versus $t_0$ for~(\ref{lin-regression}) 
		with the initial data~(\ref{initial-data}) with $\alpha = 0.1$.}
	\label{fig-2}
\end{figure}

To deal with this numerical issue, we introduce a grid of values of $t_0$ and use the linear regression~(\ref{lin-regression}) with $t_0$-dependent values of $c_1$ and the approximation error. The outcomes 
of these computations for $\alpha = 0.1$ are depicted in Figure~\ref{fig-2}, 
where the left panel shows the coefficient $c_1$ versus $t_0$ and 
the right panel shows the corresponding approximation error versus $t_0$. 
The minimal error of the size $10^{-9}$ is attained at $t_0 = 0.1738$ 
and this value of $t_0$ corresponds to $c_1 = 0.4917$, which is close 
to the claimed value $\frac{1}{2}$ in~(\ref{scaling-law}).

Using similar ideas for $u_x(t,\xi(t))$ and $u_{xx}(t,\xi(t)-0)$, we have found that the minimal approximation errors of the size $10^{-9}$ and 
$10^{-6}$ correspond to $t_0 = 0.1750$ and $t_0 = 0.1675$, respectively. 
The corresponding coefficients for the power are 
$c_1 = 1.0125$ and $c_1 = 0.4503$, which are close to the claimed values $1$ and $\frac{1}{2}$ in~(\ref{scaling-law}). It is not surprising that 
the approximation error for the second derivative $u_{xx}(t,\xi(t)-0)$ 
is significantly larger than that for the first derivative $u_x(t,\xi(t))$ 
since we use central difference approximations. Consequently, 
the coefficient $c_1 = 0.4503$ deviates from $\frac{1}{2}$ more significantly 
than the coefficient $c_1 = 1.0125$ deviates from $1$.

The accuracy is lower for larger values of $\alpha$ in the initial data~(\ref{initial-data}). For instance, computations at $\alpha = 0.5$ 
show that the linear regression~(\ref{lin-regression}) gives 
the best approximation result at $t_0 = 0.2008$ with an error of size $10^{-6}$. The coefficient $c_1 = 0.4510$ corresponds to the power $\frac{1}{2}$, which is worse than 
in the case of $\alpha = 0.1$. Similar discrepancy was found 
for $u_x(t,\xi(t))$ with the corresponding 
approximation of $c_1 = 1.0609$. 
It was surprising, however, that the accuracy of computations for 
$u_{xx}(t,\xi(t)-0)$ was comparable between the cases $\alpha = 0.1$ and $\alpha = 0.5$. The minimal error was found in the latter case of size $10^{-6}$
with corresponding coefficient $c_1 = 0.4589$.

\begin{figure}[htbp] 
	\centering
	\includegraphics[width=3in, height = 2in]{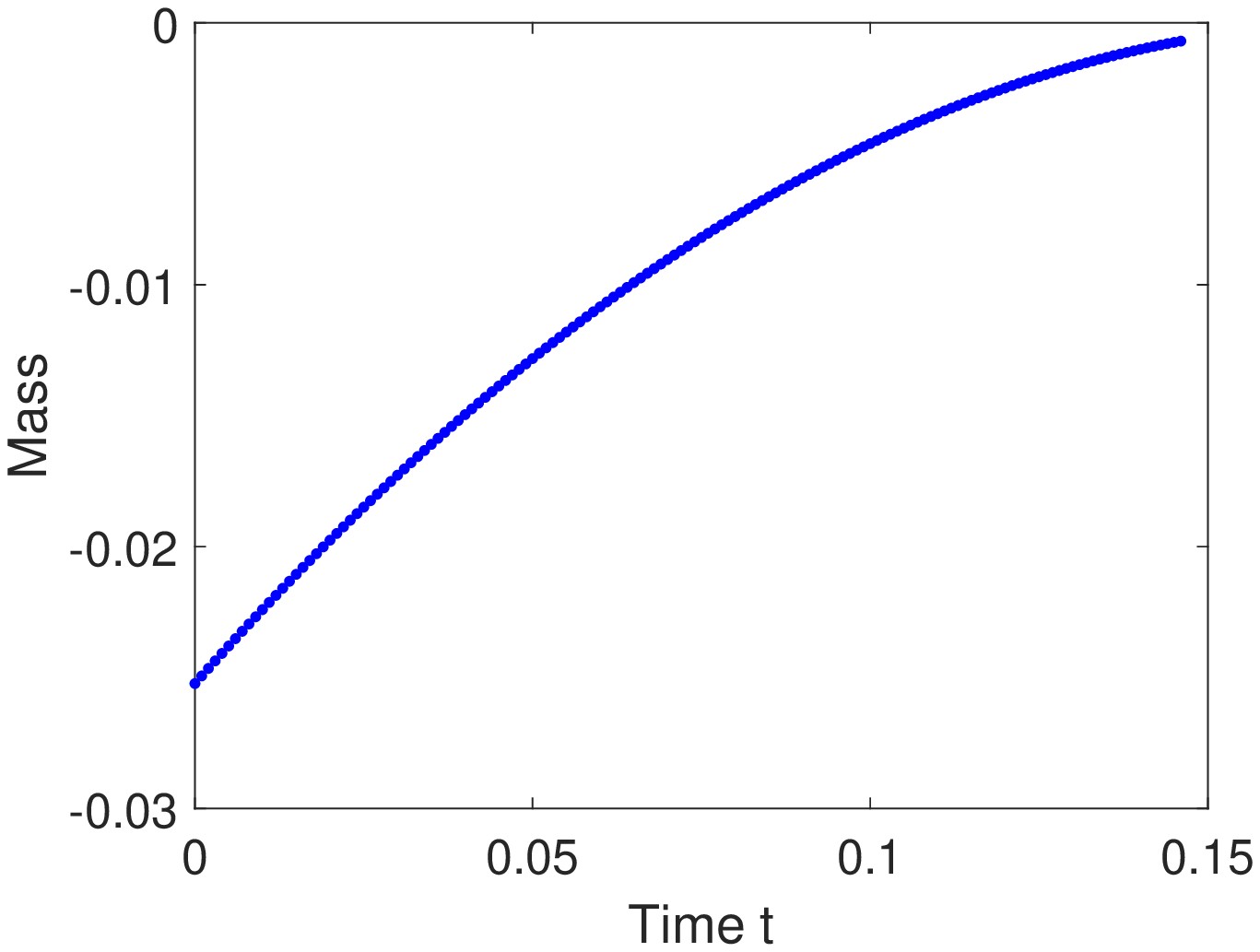}
	\includegraphics[width=3in, height = 2in]{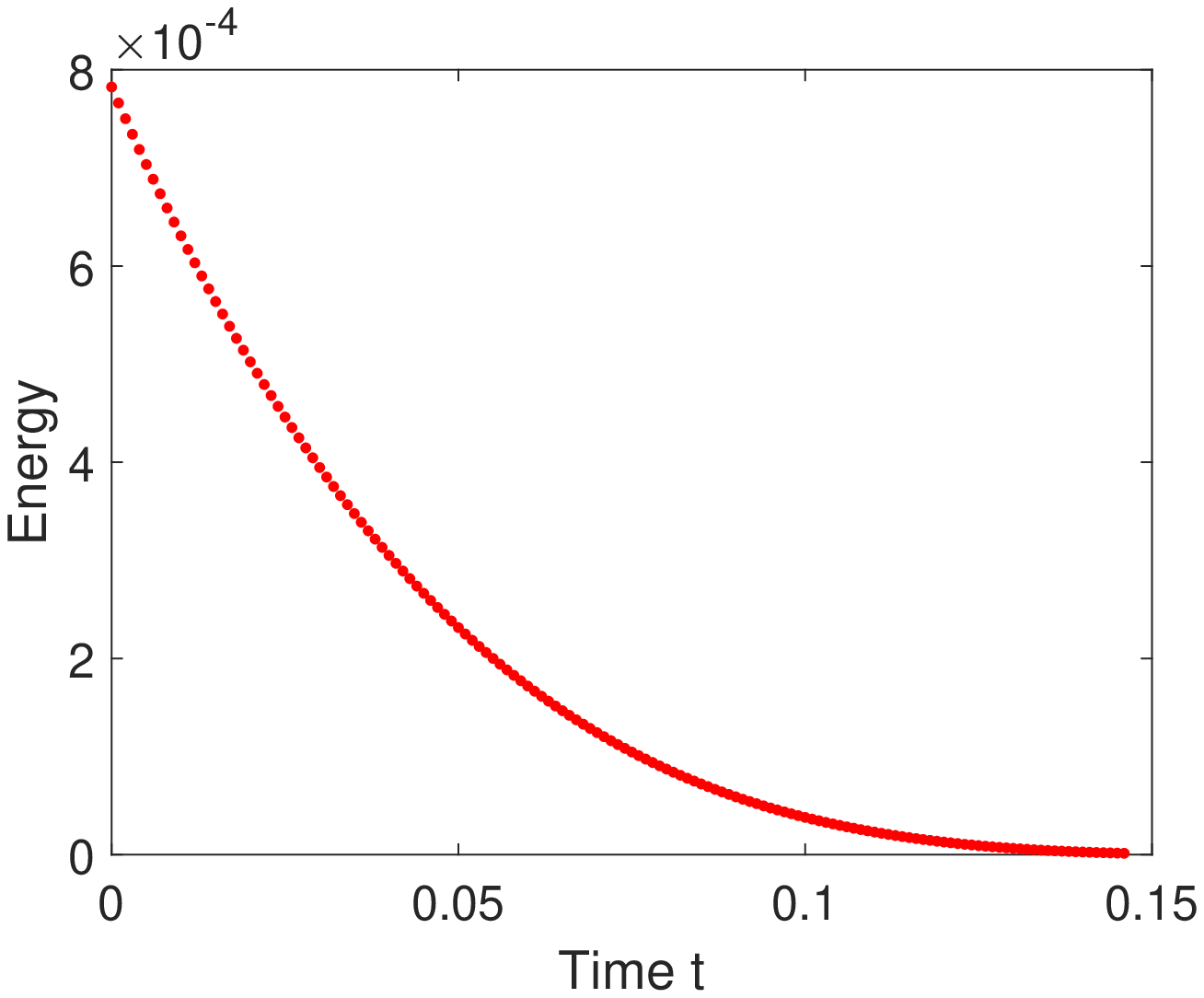} 
	\caption{Mass (left) and energy (right) versus $t$ for the time evolution 
		for the initial data~(\ref{initial-data}) with $\alpha = 0.1$.}
	\label{fig-3}
\end{figure}

We have also computed the numerical approximations for the mass and energy integrals for the compact area on $[0,\xi(t)]$, see Section~\ref{sec-3}. After the change of variables, these quantities are given by 
\begin{equation}
\mathcal M(t) := \xi(t) \int_0^1 u(t,y) dy, \quad \mathcal E(t) := \xi(t) \int_0^1 u^2(t,y) dy.
\end{equation}
Figure~\ref{fig-3} shows the evolution of the mass and energy integrals versus $t$ for the initial data~(\ref{initial-data}) with $\alpha = 0.1$. The numerically detected best power fits suggest that 
\begin{equation}
\label{scaling-law-mass}
|\mathcal M(t)| \sim (t_0-t)^2, \quad 
\mathcal E(t) \sim \sqrt{(t_0-t)^7},
\end{equation}
which are also in agreement with the balance equations~(\ref{reduction-1}) and~(\ref{reduction-2}) under the scaling laws~(\ref{scaling-law}).

\section{Conclusion}
\label{sec-6}

To summarize, we have shown analytically and numerically that the dynamics of odd viscous shocks in the modular Burgers' equations with three interfaces leads to the finite-time extinction of compact regions by means of coalescence of two consecutive interfaces. We have specified precise scaling laws for the finite-time extinction based on numerical simulations with the central difference method, which is well-adapted to deal with the nonlinear interface equations. 

These results open a road for future work to establish finite-time extinction of shocks and associated scaling laws analytically for general initial data with multiple interfaces. We anticipate that all initial data with finitely many interfaces evolve in finite time to shock waves with a single interface or to linear waves without interfaces (depending on the boundary conditions). However, it is unclear if the scaling laws, as stated in this paper, are universal for other data. 

Among other open questions, one can consider extensions of these results to the modular Burgers' equation with additional terms and to the logarithmic Burgers' equation. Dynamics of solitary waves in the modular Korteweg-de Vries equation and other Hamiltonian systems with modular nonlinearity have not been investigated so far and could also be attractive subjects of future research on their own.

\vspace{0.2cm}

{\bf Acknowledgements.} An early stage of this work was completed during the undergraduate research project of Y. Ackermann and E. Redfearn.  
The later stage of this work was completed during the visit of D. E. Pelinovsky 
to KIT as a part of Humboldt Reseach Award from Alexander von Humboldt Foundation. The project is supported by the RNF grant 19-12-00253.

\end{document}